\documentclass[11pt]{article}
\usepackage{amsfonts}
\usepackage{amssymb,amsthm,graphicx}
\usepackage{flafter}
\usepackage{epsfig}
\usepackage{cite}
\usepackage{amsmath,amscd}

\newtheorem{theorem}{Theorem}
\newtheorem{lemma}{Lemma}

\newcommand{\D}{\mbox{$\mathbb{D}$}}
\numberwithin{equation}{section}

\newtheorem{corollary}{Corollary}
\newtheorem{definition}{Definition}

\newtheorem{remark}{Remark}

\date{}

\begin{document}

\begin{center}
{\bf{\Large On third Hankel determinants for subclasses of analytic functions and close-to-convex harmonic mappings}}\\[4mm]

{\bf Yong Sun$^{1}$, Zhi-Gang Wang$^{2}$ and Antti Rasila$^{3}$}\\[2mm]

$^{1}$School of Science, Hunan Institute of Engineering,\\
Xiangtan, 411104, Hunan, People's Republic of China\\[1mm]
\textbf{E-Mail: yongsun2008@foxmail.com}\\[2mm]

$^{2}$School of Mathematics and Computing Science, Hunan First Normal University,\\
Changsha, 410205, Hunan, People's Republic of China\\[1mm]
\textbf{E-Mail: wangmath@163.com}\\[2mm]

$^{3}$Department of Mathematics and Systems Analysis, Aalto University,\\
P. O. Box 11100, FI-00076 Aalto, Finland\\[1mm]
\textbf{E-Mail: antti.rasila@iki.fi}\\[2mm]

\end{center}

\vskip.05in

\vskip.05in

\begin{center}
{\bf Abstract}
\end{center}
\begin{quotation}
In this paper, we obtain the upper bounds to the third Hankel determinants for starlike functions of order $\alpha$,
convex functions of order $\alpha$ and bounded turning functions of order $\alpha$.
Furthermore, several relevant results on a new subclass of close-to-convex harmonic mappings are obtained.
Connections of the results presented here to those that can be found in the literature are also discussed.
\end{quotation}

\vskip.02in

\noindent
{\bf 2010 Mathematics Subject Classification.} 30C45, 30C50, 58E20, 30C80.\\

\vskip.02in

\noindent
{\bf Key Words and Phrases.}
Univalent function; starlike function; convex function; bounded turning function; close-to-convex function; harmonic mapping; Hankel determinant.

\vskip.05in
\section{Introduction}
\vskip.05in

Let $\mathcal{A}$ be the class of functions {\it analytic} in the unit disk
$\D:=\{z\in\mathbb{C}:\, |z|<1\}$ of the form
\begin{equation}\label{11}
f(z)=z+\sum_{k=2}^{\infty}a_{k}z^{k}.
\end{equation}
We denote by $\mathcal{S}$ the subclass of $\mathcal{A}$ consisting of univalent functions.

\vskip.05in

A function $f\in\mathcal{A}$ is said to be starlike of order $\alpha\ (0\leq \alpha<1)$,
if it satisfies the following condition:
\begin{equation*}
\Re\bigg(\frac{zf'(z)}{f(z)}\bigg)>\alpha   \quad (z\in\D).
\end{equation*}
We denote by $\mathcal{S}^{*}(\alpha)$ the class of starlike functions of order $\alpha$.

\vskip.05in

Denote by $\mathcal{K}(\alpha)$ the class of functions $f\in\mathcal{A}$ such that
\begin{equation*}
\Re\bigg(1+\frac{zf''(z)}{f'(z)}\bigg)>\alpha   \quad (-1/2\leq\alpha<1;\ z\in\D).
\end{equation*}
In particular, functions in $\mathcal{K}(-1/2)$ are known to be close-to-convex
but are not necessarily starlike in $\D$. For $0\leq\alpha<1$,
functions in $\mathcal{K}(\alpha)$ are known to be convex of order $\alpha$ in $\D$.

\vskip.05in

A function $f\in\mathcal{A}$ is said to be in the class $\mathcal{R}(\alpha)$,
consisting of functions whose derivative have a positive real part of $\alpha\ (0\leq \alpha<1)$,
if it satisfies the following condition:
\begin{equation*}
\Re\big(f'(z)\big)>\alpha   \quad (z\in\D).
\end{equation*}

Choosing $\alpha=0$, we denote the $\mathcal{S}:=\mathcal{S}^{*}(0)$, $\mathcal{K}:=\mathcal{K}(0)$ and $\mathcal{R}:=\mathcal{R}(0)$,
the classes of starlike, convex and bounded turning functions, respectively.

\vskip.05in

Let $\mathcal{H}$ denote the class of all {\it complex-valued harmonic mappings} $f$
in $\D$ normalized by the condition $f(0)=f_{z}(0)-1=0$.
It is well-known that such functions can be written as $f=h+\overline{g}$, where $h$ and $g$ are analytic
functions in $\D$. We call $h$ the analytic part and $g$ the co-analytic part of $f$, respectively.
Let $\mathcal{S}_{H}$ be the subclass of $\mathcal{H}$ consisting of univalent and sense-preserving mappings.
Such mappings can be written in the form
\begin{equation}\label{12}
f(z)=h(z)+\overline{g(z)}=z+\sum_{k=2}^{\infty}a_{k}z^{k}
+\sum_{k=1}^{\infty}\overline{b_{k}z^{k}}
\quad (|b_{1}|<1;\ z\in\D).
\end{equation}
Harmonic mapping $f$ is called locally univalent and sense-preserving in $\D$ if and only if $|h'(z)|>|g'(z)|$ holds for $z\in\D$.
Observe that $\mathcal{S}_{H}$ reduces to $\mathcal{S}$, the class of normalized univalent analytic functions,
if the co-analytic part $g$ vanishes. The family of all functions $f\in\mathcal{S}_{H}$ with the additional property
that $f_{\overline{z}}(0)=0$ is denoted by $\mathcal{S}^{0}_{H}$.
For further information about planar harmonic mappings, see e.g. \cite{duren2004,cs1984,pr2013}.

\vskip.05in

Recall that a function $f\in\mathcal{H}$ is close-to-convex in $\D$
if it is univalent and the range $f(\D)$ is a close-to-convex domain, i.e., the
complement of $f(\D)$ can be written as the union of nonintersecting half-lines.
A normalized analytic function $f$ in $\D$ is close-to-convex in $\D$ if there
exists a convex analytic function in $\D$, not necessarily normalized, $\phi$
such that $\Re\big(f'(z)/\phi'(z)\big)>0$. In particular, if $\phi(z)=z$,
then for any $f\in\mathcal{A}$, $\Re\big(f'(z)\big)>0$ implies $f$ is close-to-convex in $\D$, see \cite {suf2005}.
We refer to \cite{bjj2013, kpv, nr1, pk1, pk2} for discussion and basic results on close-to-convex harmonic mappings.

\vskip.05in

For a harmonic mapping $f=h+\overline{g}$ in $\D$, a basic result in \cite{m2011} (see also \cite{m1981}) shows that
if at least one of the analytic functions $h$ and $g$ is convex, then $f$ is univalent whenever
it is locally univalent in $\D$. It is natural to study the univalence of $f=h+\overline{g}$ in $\D$
if it is locally univalent and sense-preserving, and analytic function $h$ is univalent and close-to-convex.
Motivated by this idea, we next consider the following subclass of $\mathcal{H}$.

\begin{definition}\label{def1}
{\rm For $\alpha\in\mathbb{C}$ with $-1/2\leq\alpha<1$, let $\mathcal{M}(\alpha)$
denote the class of harmonic mapping $f$ in $\D$ of the form \eqref{12}, with $h'(0)\neq 0$, which satisfy
\begin{equation*}
\Re\bigg(1+\frac{zh''(z)}{h'(z)}\bigg)>\alpha
\quad {\rm and} \quad  g'(z)=zh'(z) \quad \big(z\in\D\big).
\end{equation*}
}
\end{definition}

\vskip.05in

By making use of the similar arguments to those in the proof of \cite [Theorem 1]{bl2011},
one can easily obtain the close-to-convexity of the class $\mathcal{M}(\alpha)$.
For special values of $\alpha$, many authors have studied the class of close-to-convex harmonic mappings,
see e.g. \cite{m2011, chen2015, sun2016, bp, nr1}.

\vskip.05in

Pommerenke (see \cite{pommerenke1966, pommerenke1967}) defined the Hankel determinant $H_{q,n}(f)$ as
\begin{equation*}
H_{q,n}(f)=
\left|\begin{array}{cccc}
a_{n}     & a_{n+1} & \cdots & a_{n+q-1}\\
a_{n+1}   & a_{n+2} & \cdots & a_{n+q}  \\
\vdots    & \vdots  & \vdots & \vdots   \\
a_{n+q-1} & a_{n+q} & \cdots & a_{n+2(q-1)}\\
\end{array}\right|
\quad (q,n\in\mathbb{N}).
\end{equation*}

Problems involving Hankel determinants $H_{q,n}(f)$ in geometric function theory originate from the work of,
e.g., Hadamard, Polya and Edrei (see \cite{dienes,edrei}), who used them in study of singularities of meromorphic functions.
For example, they can be used in showing that a function of bounded characteristic  in $\D$, i.e.,
a function which is a ratio of two bounded analytic functions with its Laurent series around the origin having integral coefficients, is rational \cite{cantor1963}.
Pommerenke \cite{pommerenke1966} proved that the Hankel determinants of univalent functions satisfy the inequality
$|H_{q,n}(f)|<Kn^{-(\frac{1}{2}+\beta)q+\frac{3}{2}}$, where $\beta>1/4000$ and $K$ depends only on $q$.
Furthermore, Hayman \cite{hayman1968} has proved a stronger result for areally mean univalent functions, i.e.,
the estimate $H_{2,n}(f)<An^{1/2}$, where $A$ is an absolute constant.

\vskip.05in

We note that $H_{2,1}(f)$ is the well-known {\it Fekete-Szeg\H{o} functional}, see \cite{fekete1933, koepf198701, koepf198702}.
The sharp upper bounds on $H_{2,2}(f)$ were obtained by the authors of articles \cite{bansal2013, janteng2006, janteng2007, lee2013} for various classes of functions.

\vskip.05in

By the definition, $H_{3,1}(f)$ is given by
\begin{equation*}
H_{3,1}(f)=
\left|\begin{array}{ccc}
a_{1} & a_{2} & a_{3}\\
a_{2} & a_{3} & a_{4}\\
a_{3} & a_{4} & a_{5}\\
\end{array}\right|.
\end{equation*}
Note that for $f\in\mathcal{A}$, $a_{1}=1$ so that
\begin{equation*}
H_{3,1}(f)=-a_{2}^{2}a_{5}+2a_{2}a_{3}a_{4}-a_{3}^{3}+a_{3}a_{5}-a_{4}^{2}.
\end{equation*}

\vskip.05in

Obviously, the case of the upper bounds on $H_{3,1}(f)$ it is much more difficult than the cases of $H_{2,1}(f)$ and $H_{2,2}(f)$.
In 2010, Babalola \cite{babalola2010} has studied the $\max|H_{3,1}(f)|$ for the classes of starlike, convex and bounded turning functions.

\noindent{\bf Theorem A.}
\textit{ Let $f\in\mathcal{S}^{*}$, $h\in\mathcal{K}$ and $g\in\mathcal{R}$, respectively. Then
\begin{equation*}
\big|H_{3,1}(f)\big|\leq 16, \qquad
\big|H_{3,1}(h)\big|\leq \frac{32+33\sqrt{3}}{72\sqrt{3}}\approx 0.714,
\end{equation*}
and
\begin{equation*}
\big|H_{3,1}(g)\big|\leq \frac{2736\sqrt{3}+675\sqrt{5}}{4860\sqrt{3}}\approx 0.742.
\end{equation*}
}

Recently, Zaprawa \cite{zaprawa2017} proved that

\noindent{\bf Theorem B.}
\textit{ Let $f\in\mathcal{S}^{*}$, $h\in\mathcal{K}$ and $g\in\mathcal{R}$, respectively. Then
\begin{equation*}
\big|H_{3,1}(f)\big|\leq 1, \qquad
\big|H_{3,1}(h)\big|\leq \frac{49}{540}\approx 0.090,
\qquad
\big|H_{3,1}(g)\big|\leq \frac{41}{60}\approx 0.683.
\end{equation*}
}

\vskip.05in

Raza and Malik \cite{raza2013} have obtained the upper bound on $|H_{3,1}(f)|$ for a class of analytic functions that is related to the lemniscate of Bernoulli.
Also, Bansal {\it et al.} \cite{bansal2015} obtained the following results

\noindent{\bf Theorem C.}
\textit{ Let $h\in\mathcal{K}(-1/2)$ and $g\in\mathcal{R}$, respectively. Then
\begin{equation*}
\big|H_{3,1}(h)\big|\leq \frac{180+69\sqrt{15}}{32\sqrt{15}}\approx 3.609,
\qquad
\big|H_{3,1}(g)\big|\leq \frac{439}{540}\approx 0.813.
\end{equation*}
}

\vskip.05in
For the class $\mathcal{R}(\alpha)$, Vamshee Krishna {\it et al.} \cite{Vamshee2015} proved that

\noindent{\bf Theorem D.}
\textit{ Let $g\in\mathcal{R}(\alpha)$ with $0\leq \alpha\leq \frac{1}{4}$. Then
\begin{equation*}
\big|H_{3,1}(g)\big|\leq \frac{(1-\alpha)^{2}}{3}\bigg[\frac{8(1-\alpha)}{9}+\frac{1}{4}\bigg(\frac{5-4\alpha}{3}\bigg)^{\frac{3}{2}}+\frac{4}{5}\bigg].
\end{equation*}
}

\vskip.05in

In the present investigation, our goal is to discuss the upper bounds to the third Hankel determinants
for the subclasses of univalent functions: $\mathcal{S}^{*}(\alpha)$, $\mathcal{K}(\alpha)$ and $\mathcal{R}(\alpha)$.
Furthermore, we develop similar results on the Hankel determinants $|H_{3,1}(h)|$ and $|H_{3,1}(g)|$
in the context the close-to-convex harmonic mappings $f=h+\overline{g}\in\mathcal{M}(\alpha)$.

\vskip.05in
\section{Preliminary results}
\vskip.05in

Denote by $\mathcal{P}$ the class of {\it Carath\'{e}odory functions} $p$ normalized by
\begin{equation}\label{lem11}
p(z)=1+\sum_{n=1}^{\infty}p_{n}z^{n} \quad {\rm and} \quad \Re\big(p(z)\big)>0 \quad (z\in\D).
\end{equation}
Following results are the well known for functions belonging to the class $\mathcal{P}$.

\vskip.05in

\begin{lemma}\label{lem1} {\rm \cite{duren1983}}
If $p\in\mathcal{P}$ is of the form \eqref{lem11}, then
\begin{equation}\label{lem12}
|p_{n}|\leq 2 \qquad (n\in\mathbb{N}).
\end{equation}
The inequality \eqref{lem12} is sharp and the equality holds for the function
\begin{equation*}
\phi(z)=\frac{1+z}{1-z}=1+2\sum_{n=1}^{\infty}z^{n}.
\end{equation*}
\end{lemma}

\vskip.05in

\begin{lemma}\label{lem2} {\rm \cite{livingston1969}}
If $p\in\mathcal{P}$ is of the form \eqref{lem11}, then holds the sharp estimate
\begin{equation}\label{lem15}
|p_{n}-p_{k}p_{n-k}|\leq 2 \qquad (n,\ k\in\mathbb{N},\ n>k).
\end{equation}
\end{lemma}

\vskip.05in

\begin{lemma}\label{lem3} {\rm \cite{hayami2010}}
If $p\in\mathcal{P}$ is of the form \eqref{lem11}, then holds the sharp estimate
\begin{equation}\label{lem16}
|p_{n}-\mu p_{k}p_{n-k}|\leq 2 \qquad (n,\ k\in\mathbb{N},\ n>k;\ 0\leq \mu \leq 1).
\end{equation}
\end{lemma}

\vskip.05in

\begin{lemma}\label{lem4} {\rm \cite{libera1982, libera1983}}
If $p\in\mathcal{P}$ is of the form \eqref{lem11}, then there exist  $x$, $z$ such that $|x|\leq 1$ and $|z|\leq 1$,
\begin{equation}\label{lem13}
2p_{2}=p_{1}^{2}+(4-p_{1}^{2})x,
\end{equation}
and
\begin{equation}\label{lem14}
4p_{3}=p_{1}^{3}+2p_{1}(4-p_{1}^{2})x-p_{1}(4-p_{1}^{2})x^{2}+2(4-p_{1}^{2})(1-|x|^{2})z.
\end{equation}
\end{lemma}

\vskip.05in
\section{Bounds of Hankel determinants for $\mathcal{S}^{*}(\alpha)$, $\mathcal{K}(\alpha)$ and $\mathcal{R}(\alpha)$}
\vskip.05in

In this section, we assume that
\begin{equation*}
f(z)=z+\sum_{k=2}^{\infty}a_{k}z^{k}\in\mathcal{S}^{*}(\alpha),
\quad h(z)=z+\sum_{k=2}^{\infty}b_{k}z^{k}\in\mathcal{K}(\alpha),
\quad g(z)=z+\sum_{k=2}^{\infty}c_{k}z^{k}\in\mathcal{R}(\alpha).
\end{equation*}

\vskip.05in

\begin{theorem}\label{t-skr}
Let $f\in\mathcal{S}^{*}(\alpha)$, $h\in\mathcal{K}(\alpha)$ and $g\in\mathcal{R}(\alpha)$ with $0\leq \alpha<1$, respectively. Then
\begin{equation}\label{31}
\big|H_{3,1}(f)\big|\leq \frac{1}{18}(1-\alpha)^{2}(18-\alpha),
\end{equation}
\begin{equation}\label{312}
\big|H_{3,1}(h)\big|\leq \frac{1}{540}(1-\alpha)^{2}(49-16\alpha),
\end{equation}
and
\begin{equation}\label{313}
\big|H_{3,1}(g)\big|\leq \frac{1}{60}(1-\alpha)^{2}(36-20\alpha+5|1-4\alpha|).
\end{equation}
\end{theorem}

\vskip.05in

\begin{proof}
Let
\begin{equation*}
p(z)=\frac{1}{1-\alpha}\left(\frac{zf'(z)}{f(z)}-\alpha\right) \quad (0\leq \alpha<1;\ z\in\D),
\end{equation*}
then, we have $\Re\big(p(z)\big)>0$, and by elementary calculations, we obtain
\begin{equation}\label{32}
p(z)=1+\frac{1}{1-\alpha}\big(a_{2}z+(2a_{3}-a_{2}^{2})z^{2}+\cdots\big)=1+p_{1}z+p_{2}z^{2}+\ldots.
\end{equation}
It follows from \eqref{32} that
\begin{equation}\label{a2345}
\left\{\begin{array}{rcl}
 a_{2}&=&(1-\alpha)p_{1},\\
 a_{3}&=&\frac{1}{2}(1-\alpha)\big[(1-\alpha)p_{1}^{2}+p_{2}\big],\\
 a_{4}&=&\frac{1}{6}(1-\alpha)\big[(1-\alpha)^{2}p_{1}^{3}+3(1-\alpha)p_{1}p_{2}+2p_{3}\big],\\
 a_{5}&=&\frac{1}{24}(1-\alpha)\big[(1-\alpha)^{3}p_{1}^{4}+6(1-\alpha)^{2}p_{1}^{2}p_{2}+8(1-\alpha)p_{1}p_{3}+3(1-\alpha)p_{2}^{2}+6p_{4}\big].
\end{array}\right.
\end{equation}
Hence, by using the above values of $a_{2}$, $a_{3}$, $a_{4}$ and $a_{5}$ from \eqref{a2345}, and by a routine computation, we obtain
\begin{equation}\label{34}
\begin{split}
H_{3,1}(f)
&=\frac{1}{144}(1-\alpha)^{2}
\bigg\{-(1-\alpha)^{4}p_{1}^{6}+3(1-\alpha)^{3}p_{1}^{4}p_{2}+8(1-\alpha)^{2}p_{1}^{3}p_{3}-9(1-\alpha)^{2}p_{1}^{2}p_{2}^{2}\\
&\ \ \ \ -18(1-\alpha)p_{1}^{2}p_{4}+24(1-\alpha)p_{1}p_{2}p_{3}-9(1-\alpha)p_{2}^{3}+18p_{2}p_{4}-16p_{3}^{2}\bigg\}.
\end{split}
\end{equation}
From \eqref{34}, we have
\begin{equation*}
\begin{split}
H_{3,1}(f)
&=\frac{1}{144}(1-\alpha)^{2}
\bigg\{(1-\alpha)\big[p_{2}-(1-\alpha)p_{1}^{2}\big]^{3}-16\big[p_{3}-(1-\alpha)p_{1}p_{2}\big]^{2}\\
&\ \ \ \ +8\big[p_{2}-(1-\alpha)p_{1}^{2}\big]\big[p_{4}-(1-\alpha)p_{1}p_{3}\big]
+10\big[p_{2}-(1-\alpha)p_{1}^{2}\big]\big[p_{4}-(1-\alpha)p_{2}^{2}\big]\bigg\}.
\end{split}
\end{equation*}
We note that
\begin{equation*}
0<1-\alpha\leq 1  \quad {\rm for} \quad  0\leq \alpha<1,
\end{equation*}
by triangle inequality and Lemma \ref{lem3}, we obtain the estimate \eqref{31} of $H_{3,1}(f)$.

Next, we consider $H_{3,1}(h)$. According to the Alexander relation, $b_{k}=ka_{k}\ (k\in\mathbb{N})$.
Putting it into the definition of $H_{3,1}(h)$ and applying the formula \eqref{a2345}, we have
\begin{equation}\label{k1}
\begin{split}
H_{3,1}(h)
&=\frac{1}{8640}(1-\alpha)^{2}
\bigg\{-(1-\alpha)^{4}p_{1}^{6}+6(1-\alpha)^{3}p_{1}^{4}p_{2}+12(1-\alpha)^{2}p_{1}^{3}p_{3}-21(1-\alpha)^{2}p_{1}^{2}p_{2}^{2}\\
&\ \ \ \ -36(1-\alpha)p_{1}^{2}p_{4}+36(1-\alpha)p_{1}p_{2}p_{3}-4(1-\alpha)p_{2}^{3}+72p_{2}p_{4}-60p_{3}^{2}\bigg\}.
\end{split}
\end{equation}
From \eqref{k1}, we have
\begin{equation*}
\begin{split}
H_{3,1}(h)
&=\frac{1}{8640}(1-\alpha)^{2}
\bigg\{8(1-\alpha)\big[p_{2}-\frac{1}{2}(1-\alpha)p_{1}^{2}\big]^{3}+24p_{4}\big[p_{2}-(1-\alpha)p_{1}^{2}\big]\\
&\ \ \ \ +36p_{2}\big[p_{4}-(1-\alpha)p_{2}^{2}\big]+12\big[p_{2}-(1-\alpha)p_{1}^{2}\big]\big[p_{4}-(1-\alpha)p_{1}p_{3}\big]\\
&\ \ \ \ -60p_{3}\big[p_{3}-\frac{4}{5}(1-\alpha)p_{1}p_{2}\big]+24(1-\alpha)p_{2}^{2}\big[p_{2}-\frac{3}{8}(1-\alpha)p_{1}^{2}\big]\bigg\}.
\end{split}
\end{equation*}
We observe that for $0\leq\alpha<1$ holds
\begin{equation*}
1-\alpha, \quad \frac{1}{2}(1-\alpha), \quad \frac{4}{5}(1-\alpha),\quad \frac{3}{8}(1-\alpha)\in [0,1].
\end{equation*}
By using Lemma \ref{lem1} and Lemma \ref{lem3} and triangle inequality, it easy to get the estimate \eqref{312} of $H_{3,1}(g)$.

Finally, for $H_{3,1}(g)$. Let
\begin{equation*}
\frac{1}{1-\alpha}\left(g'(z)-\alpha\right)=1+\sum_{k=1}^{\infty}p_{k}z^{k}\in\mathcal{P}.
\end{equation*}
If $g\in\mathcal{R}(\alpha)$, then
\begin{equation}\label{cn}
(k+1)c_{k+1}=(1-\alpha)p_{k}  \qquad      (k\in\mathbb{N}).
\end{equation}
Putting it into the definition of $H_{3,1}(g)$ and by the same way, we have
\begin{equation*}
\begin{split}
H_{3,1}(g)
&=\frac{1}{2160}(1-\alpha)^{2}
\bigg\{(1-\alpha)\big[-108p_{1}^{2}p_{4}+180p_{1}p_{2}p_{3}-80p_{2}^{3}\big]+144p_{2}p_{4}-135p_{3}^{2}\bigg\}\\
&=\frac{1}{2160}(1-\alpha)^{2}
\bigg\{108(1-\alpha)p_{4}(p_{2}-p_{1}^{2})+80(1-\alpha)p_{2}(p_{4}-p_{2}^{2})\\
&\ \ \ \ -135p_{3}(p_{3}-p_{1}p_{2})-45(1-4\alpha)p_{2}(p_{4}-p_{1}p_{3})+(1+8\alpha)p_{2}p_{4}\bigg\}.
\end{split}
\end{equation*}
Hence, it is easy to obtain the bound of $H_{3,1}(g)$. This completes the proof.
\end{proof}

\vskip.05in

\begin{remark}
{\rm By setting $\alpha=0$ in Theorem \ref{t-skr}, we obtain the known results of Theorem { B}, and they are much better than Theorem { A}.
Furthermore, the bounds of $H_{3,1}(g)$ in \eqref{313} improved and extended the result of the Theorem { D}.
}
\end{remark}

\vskip.05in

In 1960, Lawrence Zalcman posed a conjecture that the coefficients of $\mathcal{S}$ satisfy the sharp inequality
\begin{equation*}
|a_{n}^{2}-a_{2n-1}| \leq (n-1)^{2}  \qquad (n\in\mathbb{N}),
\end{equation*}
with equality only for the Koebe function $k(z)=z/(1-z)^{2}$ and its rotations.
We call $J_{n}(f)=a_{n}^{2}-a_{2n-1}$ the Zalcman functional for $f\in\mathcal{S}$.

We observe that $H_{3,1}(f)\ (f\in\mathcal{A})$ can be written in the form
\begin{equation*}
H_{3,1}(f)=a_{3}(a_{2}a_{4}-a_{3}^{2})+a_{4}(a_{2}a_{3}-a_{4})-a_{5}J_{2}(f),
\end{equation*}
and equivalently,
\begin{equation*}
H_{3,1}(f)=a_{3}J_{3}(f)+a_{4}(2a_{2}a_{3}-a_{4})-a_{5}a_{2}^{2}.
\end{equation*}

\vskip.05in

An analogous calculation can be applied to the Zalcman functional $J_{n}(f)$ for the classes of starlike, convex and bounded turning functions of order $\alpha$.

\vskip.05in

\begin{theorem}\label{t5-32} The following estimates hold for $J_{n}(f)$:
\begin{enumerate}
\item
If $f\in\mathcal{S}^{*}(\alpha)\ (0\leq \alpha<1)$, then $J_{3}(f)\leq \frac{1}{2}(1-\alpha)(8-7\alpha)$.
\item
If $h\in\mathcal{K}(\alpha)\ (-1/2\leq \alpha<1)$, then $J_{3}(h)\leq \frac{1}{360}(1-\alpha)(127-109\alpha)$.
\item
If $g\in\mathcal{R}(\alpha)\ (0\leq \alpha<1)$, then $J_{n}(g)\leq \frac{2}{2n-1}(1-\alpha)  \qquad (n\geq 2)$.
\end{enumerate}

\end{theorem}

\vskip.05in

\begin{proof}
Let $f\in\mathcal{S}^{*}(\alpha)$, from \eqref{a2345}, it follow that
\begin{equation*}
\begin{split}
J_{3}(f)&=\frac{1}{24}(1-\alpha)
\bigg\{-5(1-\alpha)^{3}p_{1}^{4}-6(1-\alpha)^{2}p_{1}^{2}p_{2}-3(1-\alpha)p_{2}^{2}+8(1-\alpha)p_{1}p_{3}+6p_{4}\bigg\}\\
&=\frac{1}{24}(1-\alpha)
\bigg\{-5(1-\alpha)\big[p_{2}-(1-\alpha)p_{1}^{2}\big]^{2}+8(1-\alpha)p_{1}\big[p_{3}-(1-\alpha)p_{1}p_{2}\big]\\
&\ \ \ \ +8(1-\alpha)p_{2}\big[p_{2}-(1-\alpha)p_{1}^{2}\big]+6\big[p_{4}-(1-\alpha)p_{2}^{2}\big]\bigg\}.
\end{split}
\end{equation*}
By using Lemma \ref{lem1} and Lemma \ref{lem3}, we obtain the above bound for the Zalcman functional $J_{3}(f)$.

\vskip.05in

Combining the Alexander relation, $b_{k}=ka_{k}$, and the formula \eqref{a2345}, yields
\begin{equation*}
\begin{split}
J_{3}(h)&=\frac{1}{360}(1-\alpha)
\bigg\{-7(1-\alpha)^{3}p_{1}^{4}-2(1-\alpha)^{2}p_{1}^{2}p_{2}-(1-\alpha)p_{2}^{2}+24(1-\alpha)p_{1}p_{3}+18p_{4}\bigg\}\\
&=\frac{1}{360}(1-\alpha)
\bigg\{-\frac{63}{4}(1-\alpha)\big[p_{2}-\frac{2}{3}(1-\alpha)p_{1}^{2}\big]^{2}+24(1-\alpha)p_{1}\big[p_{3}-\frac{2}{3}(1-\alpha)p_{1}p_{2}\big]\\
&\ \ \ \ +\frac{21}{2}(1-\alpha)p_{2}\big[p_{2}-\frac{2}{3}(1-\alpha)p_{1}^{2}\big]+\frac{17}{4}(1-\alpha)p_{2}^{2}+18p_{4}\bigg\}.
\end{split}
\end{equation*}
Again, by using Lemma \ref{lem1} and Lemma \ref{lem3}, we obtain the bound for the Zalcman functional $J_{3}(h)$.

\vskip.05in

For $g\in\mathcal{R}(\alpha)$, according to the formula \eqref{cn}, we have
\begin{equation*}
\begin{split}
J_{n}(g)&=\frac{1}{n^{2}}(1-\alpha)^{2}p_{n-1}^{2}-\frac{1}{2n-1}(1-\alpha)p_{2n-2}\\
&=-\frac{1}{2n-1}(1-\alpha)\bigg[p_{2n-2}-\frac{2n-1}{n^{2}}(1-\alpha)p_{n-1}^{2}\bigg].
\end{split}
\end{equation*}
In view of
\begin{equation*}
0<\frac{2n-1}{n^{2}}(1-\alpha)<1    \qquad (0\leq \alpha<1;\ n\geq 2),
\end{equation*}
and, by Lemma \ref{lem3}, we have the desired bound of the Zalcman functional $J_{n}(g)$. This completes the proof.
\end{proof}

\vskip.05in

\begin{remark}
{\rm By setting $\alpha=-1/2$ for the class $\mathcal{K}(\alpha)$ in Theorem \ref{t5-32}, we obtain the known results \cite [Theorem 2.3]{abu2014}.
Furthermore, using the similar argument in Theorem \ref{t5-32}, we may obtain the bounds of the Zalcman functional $J_{2}(f)$ and $J_{2}(h)$:
If $f\in\mathcal{S}^{*}(\alpha)\ (0\leq \alpha<1)$, then $J_{2}(f)\leq 1-\alpha$. If
$h\in\mathcal{K}(\alpha)\ (-1/2\leq \alpha<1)$, then $J_{2}(h)\leq \frac{1}{3}(1-\alpha)$ .
}
\end{remark}

\vskip.05in

\vskip.05in
\section{Bounds of Hankel determinants for $\mathcal{M}(\alpha)$}
\vskip.05in

In this section, we obtain upper bounds for the Hankel determinants $|H_{3,1}(h)|$ and $|H_{3,1}(g)|$
of close-to-convex harmonic mappings $f=h+\overline{g}\in\mathcal{M}(\alpha)$.

\vskip.05in

\begin{theorem}\label{tH31hg}
Let $f=h+\overline{g}\in\mathcal{M}(\alpha)$ be of the form \eqref{12}. Then
\begin{equation*}
\big|H_{3,1}(h)\big|
\leq \frac{1}{540}(1-\alpha)^{2}(15\alpha^{2}-34\alpha+52),
\end{equation*}
and
\begin{equation*}
\big|H_{3,1}(g)\big|\leq \frac{1}{30}(1-\alpha).
\end{equation*}

\end{theorem}

\vskip.05in

\begin{proof}
Let
\begin{equation*}
p(z)=\frac{1}{1-\alpha}\left(1+\frac{zh''(z)}{h'(z)}-\alpha\right)
=1+\sum_{k=1}^{\infty}p_{k}z^{k}\in\mathcal{P} \qquad (-\frac{1}{2}\leq \alpha<1;\ z\in\D).
\end{equation*}

Using the same method of Theorem \ref{t-skr}, we get the expression of $H_{3,1}(h)$ is the formula \eqref{k1}.
We give another decomposition for functional $H_{3,1}(h)$ as follows
\begin{equation*}
\begin{split}
H_{3,1}(h)
&=\frac{1}{8640}(1-\alpha)^{2}
\bigg\{8(1-\alpha)\big[p_{2}-\frac{1}{2}(1-\alpha)p_{1}^{2}\big]^{3}-60\big[p_{3}-\frac{1}{2}(1-\alpha)p_{1}p_{2}\big]^{2}\\
&\ \ \ \ +48\big[p_{2}-\frac{1}{2}(1-\alpha)p_{1}^{2}\big]\big[p_{4}-\frac{1}{2}(1-\alpha)p_{1}p_{3}\big]
-15(1-\alpha)^{2}p_{1}^{2}p_{2}^{2}\\
&\ \ \ \ +24\big[p_{2}-\frac{1}{2}(1-\alpha)p_{1}^{2}\big]\big[p_{4}-\frac{1}{2}(1-\alpha)p_{2}^{2}\big]\bigg\}.
\end{split}
\end{equation*}
We note that
\begin{equation*}
0\leq \frac{1}{2}(1-\alpha)\leq 1  \quad {\rm for} \quad  0\leq \alpha<1,
\end{equation*}
by triangle inequality and Lemmas \ref{lem1}-\ref{lem3}, we can obtain the estimate of $H_{3,1}(h)$.

\vskip.05in

By the power series representations of $h$ and $g$ for $f=h+\overline{g}\in\mathcal{M}(\alpha)$, we see that
\begin{equation*}
b_{1}=0, \qquad (k+1)b_{k+1}=k a_{k}  \quad {\rm for} \quad k\geq 1,
\end{equation*}
which yields
\begin{equation*}
\left\{\begin{array}{rcccl}
 b_{2}&=&\frac{1}{2}a_{1}&=&\frac{1}{2}, \\
 b_{3}&=&\frac{2}{3}a_{2}&=&\frac{1}{3}(1-\alpha)p_{1}, \\
 b_{4}&=&\frac{3}{4}a_{3}&=&\frac{1}{8}\big[(1-\alpha)^{2}p_{1}^{2}+(1-\alpha)p_{2}\big],\\
 b_{5}&=&\frac{4}{5}a_{4}&=&\frac{1}{30}\big[(1-\alpha)^{3}p_{1}^{3}+3(1-\alpha)^{2}p_{1}p_{2}+2(1-\alpha)p_{3}\big].
\end{array}\right.
\end{equation*}
Then, by using \eqref{lem13} and \eqref{lem14} in Lemma \ref{lem4}, we obtain that for some $x$ and $z$ such that $|x|\leq 1$ and $|z|\leq 1$,
\begin{equation*}
\begin{split}
H_{3,1}(g)
& =2b_{2}b_{3}b_{4}-b_{3}^{3}-b_{2}^{2}b_{5}=b_{3}b_{4}-b_{3}^{3}-\frac{1}{4}b_{5}\\
& =\frac{1}{2160}(1-\alpha)\bigg\{\big(-8\alpha^{2}+16\alpha+1\big)p_{1}^{3}+9(4-p_{1}^{2})\big[p_{1}x^{2}-2(1-|x|^{2})z\big]\bigg\}.
\end{split}
\end{equation*}
By Lemma \ref{lem1}, we may assume that $|p_{1}|=c\in[0,2]$.
By applying the triangle inequality in above relation with $\mu=|x|$, we obtain
\begin{equation*}
\big|H_{3,1}(g)\big|
\leq \frac{1}{2160}(1-\alpha)\bigg\{\big|8\alpha^{2}-16\alpha-1\big|c^{3}+9(4-c^{2})\big[(c-2)\mu^{2}+2\big]\bigg\}
=:Q(c,\mu).
\end{equation*}
We note that
\begin{equation*}
(c-2)\mu^{2}+2\leq 2 , \quad {\rm for} \quad \mu\in[0,1] \quad {\rm and} \quad  c\in[0,2].
\end{equation*}
Hence, we have
\begin{equation*}
\big|H_{3,1}(g)\big|\leq Q(c,\mu)\leq Q(c,0)=\frac{1}{2160}(1-\alpha)\bigg\{\big|8\alpha^{2}-16\alpha-1\big|c^{3}-18c^{2}+72\bigg\}.
\end{equation*}
Let
\begin{equation*}
\chi(c)=\big|8\alpha^{2}-16\alpha-1\big|c^{3}-18c^{2}+72 \quad (c\in[0, 2]).
\end{equation*}
Then, we obtain
\begin{equation*}
\chi'(c)=3c\big(\big|8\alpha^{2}-16\alpha-1\big|c-12\big),
\end{equation*}
and
\begin{equation*}
\chi''(c)=6\big(\big|8\alpha^{2}-16\alpha-1\big|c-6\big).
\end{equation*}
Solving the equation $\chi'(c)=0$, we get the critical points are $c=0$ and
\begin{equation*}
C_{1}=\frac{12}{\big|8\alpha^{2}-16\alpha-1\big|}.
\end{equation*}
We observe that
\begin{equation*}
\chi''(c)\Big|_{c=0}=-36<0, \qquad \chi''(c)\Big|_{c=C_{1}}=36>0,
\end{equation*}
and
\begin{equation*}
0\leq\big|8\alpha^{2}-16\alpha-1\big|\leq 9 \quad (-1/2\leq \alpha<1).
\end{equation*}
Hence, we get
\begin{equation*}
\chi(c)\leq \max\bigg\{\chi(0),\ \chi(2)\bigg\}=\max\bigg\{72,\ 8\big|8\alpha^{2}-16\alpha-1\big|\bigg\}=72.
\end{equation*}
Thus, we obtain the following bound
\begin{equation*}
\big|H_{3,1}(g)\big|\leq \frac{1}{30}(1-\alpha).
\end{equation*}

\end{proof}

\vskip.05in

\begin{remark}
{\rm In order to obtain the bounds of $H_{3,1}(h)$, we give two kinds of decomposition for formula (3.7)
in Theorem \ref{t-skr} and Theorem \ref{tH31hg}, respectively. Hence, it is a natural question:
Whether there is an optimal decomposition for the similar formulae.
}
\end{remark}

\vskip.05in

\begin{remark}
{\rm For $H_{3,1}(g)$ in Theorem \ref{tH31hg}, if we apply the method in Theorem \ref{t-skr}, then
\begin{equation*}
\begin{split}
H_{3,1}(g)
&=2b_{2}b_{3}b_{4}-b_{3}^{3}-b_{2}^{2}b_{5}=b_{3}b_{4}-b_{3}^{3}-\frac{1}{4}b_{5}\\
&=\frac{1}{540}(1-\alpha)\bigg\{-2(1-\alpha)^{2}p_{1}^{3}-9\big[p_{3}-(1-\alpha)p_{1}p_{2}\big]\bigg\}\\
&=\frac{1}{540}(1-\alpha)\bigg\{3(1-\alpha)p_{1}\big[p_{2}-\frac{2}{3}(1-\alpha)p_{1}^{2}\big]-9\big[p_{3}-\frac{2}{3}(1-\alpha)p_{1}p_{2}\big]\bigg\}.
\end{split}
\end{equation*}
By using Lemmas \ref{lem1} and \ref{lem3}, we have
\begin{equation*}
\big|H_{3,1}(g)\big|\leq \frac{1}{90}(1-\alpha)(5-2\alpha),
\end{equation*}
obviously,
\begin{equation*}
\frac{1}{90}(1-\alpha)(5-2\alpha)> \frac{1}{30}(1-\alpha)  \quad {\rm for} \quad -\frac{1}{2}\leq \alpha<1.
\end{equation*}
Hence, we choose the bound of in $H_{3,1}(g)$ in Theorem \ref{tH31hg}.
}
\end{remark}

\vskip.05in

\begin{corollary}\label{cor331}
Let $f=h+\overline{g}\in\mathcal{M}(-1/2)$ be of the form \eqref{12}. Then
\begin{equation*}
\big|H_{3,1}(h)\big|\leq \frac{291}{960}\approx 0.303125,   \qquad   \big|H_{3,1}(g)\big|\leq \frac{1}{20}=0.05.
\end{equation*}
\end{corollary}

\vskip.05in

\begin{remark}
{\rm The result of $H_{3,1}(h)$ in Corollary \ref{cor331} is much better than Theorem {\bf C} (see \cite [Theorem 2.7]{bansal2015}).
From the upper bounds of $H_{3,1}(h)$ and $H_{3,1}(g)$, we note that the former is much larger than the latter,
this implies that the analytic part $h$ accounts for absolute advantage than the co-analytic part $g$
for the harmonic mappings $f=h+\overline{g}\in\mathcal{M}(\alpha)$.
}
\end{remark}

\vskip.10in
\noindent {\large\bf{Acknowledgements}}
\vskip.05in

The present investigation was supported by the \textit{Natural
Science Foundation of Hunan Province} under Grant no. 2016JJ2036, and the \textit{Foundation of Educational Committee of Hunan Province} under Grant no. 15C1089.

\end{document}